\documentclass[12pt]{amsart}

\usepackage{mathtools,amsmath,amsthm,amsfonts,latexsym,amssymb,mathrsfs,setspace,verbatim,cite,mathtools,tikz-cd,appendix,comment,enumitem,url}
\usepackage{fullpage}

\mathtoolsset{showonlyrefs=true}

\newtheorem{thm}{Theorem}

\newtheorem{lem}[thm]{Lemma}

\newtheorem{cor}[thm]{Corollary}
\numberwithin{thm}{section}
\numberwithin{equation}{section}

\theoremstyle{definition}

\newtheorem*{ex}{Example}

\newcommand{\rat}{\mathbb Q}
\newcommand{\real}{\mathbb R}

\newcommand{\alg}{\overline\rat}
\newcommand{\algt}{\alg^{\times}}

\newcommand{\intg}{\mathbb Z}
\newcommand{\nat}{\mathbb N}

\newcommand{\tors}{\mathrm{tors}}
\newcommand{\divv}{\mathrm{div}}

\newcommand{\spann}{\mathrm{span}}
\newcommand{\norm}{\mathrm{Norm}}

\newcommand{\rank}{\mathrm{rank}}
\newcommand{\bb}{{\bf b}}
\newcommand{\yy}{{\bf y}}
\newcommand{\xx}{{\bf x}}

\begin{document}

%\title[Extensions of completely additive functions]{Extensions of completely additive arithmetic functions to $\overline{\mathbb Q}^\times/\overline{\mathbb Q}^\times_{\mathrm{tors}}$}

\title[A classification of linear maps]{A classification of $\mathbb Q$-linear maps from $\overline{\mathbb Q}^\times/\overline{\mathbb Q}^\times_{\mathrm{tors}}$ to $\mathbb R$}

\author[C.L. Samuels]{Charles L. Samuels}
\address{Christopher Newport University, Department of Mathematics, 1 Avenue of the Arts, Newport News, VA 23606, USA}
\email{charles.samuels@cnu.edu}
\thanks{This work was funded in part by the AMS-Simons Research Enhancement Grant for PUI Faculty}

\begin{abstract}
	A 2009 article of Allcock and Vaaler explored the $\mathbb Q$-vector space $\mathcal G := \overline{\mathbb Q}^\times/{\overline{\mathbb Q}^\times_{\mathrm{tors}}}$,
	showing how to represent it as part of a function space on the places of $\overline{\mathbb Q}$.  We establish a representation theorem for the $\mathbb R$-vector space 
	of $\mathbb Q$-linear maps from $\mathcal G$ to $\mathbb R$, enabling us to classify extensions to $\mathcal G$ of completely additive arithmetic functions.
	We further outline a strategy to construct $\mathbb Q$-linear maps from $\mathcal G$ to $\mathbb Q$, i.e., elements of the algebraic dual of $\mathcal G$.  
	Our results make heavy use of Dirichlet's $S$-unit Theorem as well as a measure-like object called a consistent map, first introduced by the author in previous work.
\end{abstract}

\maketitle

\section{Introduction}

\subsection{Background} Let $\alg$ be a fixed algebraic closure of $\rat$ and let $\algt_\tors$ denote the group of roots of unity in $\algt$.  We write $\mathcal G = \algt/\algt_\tors$
and note that $\mathcal G$ is a vector space over $\rat$ with addition and scalar multiplication given by
\begin{equation} \label{GOperations}
	(\alpha,\beta) \mapsto \alpha\beta\quad\mbox{and}\quad (r,\alpha) \mapsto \alpha^r.
\end{equation}
An innovative article of Allcock and Vaaler \cite{AllcockVaaler} showed how to interpret $\mathcal G$ as a certain function space in the following way.

For each number field $K$, we write $M_K$ to denote the set of all places of $K$.   If $L/K$ is a finite extension and $w\in M_L$, then $w$ divides a unique place $v$ of $K$, and in this case, 
we write $K_w$ to denote the completion of $K$ with respect to $v$.  Additionally, we let $p_v$ be the unique place of $\rat$ such that $v$ divides $p_v$ and
let $\|\cdot \|_v$ be the unique extension to $K_v$ of the usual $p_v$-adic absolute value on $\rat_v$.  In this notation, the well-known product formula may be expressed as
\begin{equation*}
	\sum_{v\in M_K} \frac{[K_v:\rat_v]}{[K:\rat]} \log \|\alpha\|_v = 0
\end{equation*}
for all non-zero elements $\alpha\in K$.

Letting $Y$ denote the set of all places of $\alg$, we define $Y(K,v) = \{y\in Y:y\mid v\}$.  Further setting
\begin{equation} \label{Jdef}
	\mathcal J = \left\{(K,v):[K:\rat]<\infty,\ v\in M_K\right\},
\end{equation}
Allcock and Vaaler observed that the collection $\{Y(K,v): (K,v)\in \mathcal J\}$ is a basis for a totally disconnected, Hausdorff topology on $Y$, and moreover, there is a Borel measure $\lambda$ on $Y$ such that
\begin{equation*}
	\lambda(Y(K,v)) = \frac{[K_v:\rat_v]}{[K:\rat]}\quad\mbox{for all } (K,v)\in \mathcal J.
\end{equation*}
Each element $\alpha\in \mathcal G$ corresponds to a locally constant function $f_\alpha:Y\to \real$ given by the formula $f_\alpha(y) = \log \|\alpha\|_y$.  When $\mathcal G$ is equipped with a
norm arising from the Weil height, they proved that $\alpha\mapsto f_\alpha$ defines an isometric isomorphism from $\mathcal G$ onto a dense $\rat$-linear subspace of
\begin{equation*}
	\left\{f\in L^1(Y): \int_Y f(y)d\lambda(y) = 0\right\}.
\end{equation*}

More recently, the author \cite{SamuelsClassification,SamuelsConsistent} began the study of various dual spaces related to $\mathcal G$.  To that end, we defined a map $c:\mathcal J\to \real$ to be {\it consistent} if 
\begin{equation} \label{Consistency}
	c(K,v) = \sum_{w\mid v} c(L,w)
\end{equation}
for all number fields $K$, all places $v$ of $K$, and all finite extensions $L/K$.  The set of all consistent maps forms a vector space over $\real$ with addition and scalar multiplication given by the formulas
\begin{equation*}
	(c+d)(K,v) = c(K,v) + d(K,v) \quad\mbox{and}\quad (rc)(K,v) = rc(K,v).
\end{equation*}
We shall write $\mathcal J^*$ to denote this space.  Every Radon measure $\mu$ on $Y$ yields a corresponding consistent map via the formula $c(K,v):=\mu(Y(K,v))$, however, not every consistent map is built in this way
(see \cite{AbergSamuels}).  The most fundamental consistent map arises from the measure $\lambda$ appearing in \cite{AllcockVaaler} which we shall simply denote by $\lambda$, ie.,
\begin{equation*}
	\lambda(K,v) = \frac{[K_v:\rat_v]}{[K:\rat]}.
\end{equation*}
The main result of \cite{SamuelsClassification} constructed an isomorphism between $\algt/\overline \intg^\times$ and a certain $\rat$-linear subspace of $\mathcal J^*$.
Later in \cite{SamuelsConsistent}, we studied the $\real$-vector space $LC_c(Y)$ of locally constant functions from $Y$ to $\real$ with compact support, and additionally, we examined its co-dimension $1$ subspace
\begin{equation*}
	LC_c^0(Y) = \left\{f\in LC_c(Y): \int_Yf(y)d\lambda(y) = 0\right\}.
\end{equation*}
As a special case of a more general set of theorems, we showed that
\begin{equation} \label{LCIsos}
	\mathcal J^*\cong LC_c(Y)^* \quad\mbox{and}\quad \mathcal J^*/\spann_\real\{\lambda\} \cong LC_c^0(Y)^*.
\end{equation}
In both cases, these isomorphisms are defined explicitly, and as such, we regard them as algebraic versions of the Riesz Representation Theorem (see \cite{Rudin,DunfordSchwartz,Bogachev}, for example). 

\subsection{Main Results} While we consider \cite{SamuelsClassification,SamuelsConsistent} to be strong results, they leave open any questions about two important spaces:
\begin{enumerate}[label=(\Roman*)]
	\item\label{AlgDual} The $\rat$-vector space of $\rat$-linear maps from $\mathcal G$ to $\rat$, i.e., the algebraic dual $\mathcal G^*$ of $\mathcal G$.
	\item\label{RealMaps} The $\real$-vector space of $\rat$-linear maps from $\mathcal G$ to $\real$, which we shall denote by $\mathcal L(\mathcal G, \real)$.
\end{enumerate}
For each $c\in \mathcal J^*$ and $\alpha\in \alg^\times$, we let
\begin{equation*}
	\Phi_c(\alpha) = \sum_{v\in M_K} c(K,v)\log\|\alpha\|_v,
\end{equation*}
where $K$ is any number field containing $\alpha$.  Because $c$ is assumed to be consistent, this definition does not depend on the choice of $K$, and moreover, its value is unchanged when $\alpha$ is multiplied by a root of unity.
Hence, $\Phi_c:\mathcal G \to \real$ is a well-defined $\rat$-linear map, and we may define $\Phi^*:\mathcal J^*\to \mathcal L(\mathcal G,\real)$ by $\Phi^*(c) = \Phi_c$.  Our main result is the following representation theorem for 
$\mathcal L(\mathcal G, \real)$.  

\begin{thm} \label{RealSurjective}
	The map $\Phi^*:\mathcal J^* \to \mathcal L(\mathcal G,\real)$ is a surjective $\real$-linear transformation such that $\ker(\Phi^*) = \spann_\real\{\lambda\}$.
\end{thm}

If $F$ is a number field and $q$ is a place of $F$, we observe that $\mathcal J_q^*:= \left\{ c\in \mathcal J^*: c(F,q) = 0\right\}$ is a subspace of $\mathcal J^*$.
Given an arbitrary consistent map $c\in \mathcal J^*$, the coset $c + \ker(\Phi^*)$ contains a unique element $d\in \mathcal J_q^*$, namely
\begin{equation*}
	d(K,v) = c(K,v) - c(F,q)\lambda(K,v).
\end{equation*}
This observation yields the following consequence of Theorem \ref{RealSurjective}.

\begin{cor} \label{RestrictedIso}
	Let $F$ be a number field and let $q$ be a place of $F$.   Then the map $c\mapsto \Phi_c$ defines an $\real$-vector space isomorphism from $\mathcal J_q^*$ to $\mathcal L(\mathcal G,\real)$.
\end{cor}

By using the specific case $F=\rat$ and $q=\infty$, Corollary \ref{RestrictedIso} provides a useful framework to classify extensions to $\mathcal G$ of completely additive arithmetic functions.  
Assuming that $c\in \mathcal J_\infty^*$, we note the following famous examples:
\begin{description}%[label=(AF.\arabic*)]
	\item[Natural Logarithm] $\Phi_c$ extends the natural logarithm on $\nat$ if and only if $c(\rat,p) = -1$ for all $p\ne \infty$.
	\item[Prime Omega Function] Let $\Omega(n)$ be the number of prime factors of $n$, counted with multiplicity (see \cite{Lalin,Dressler,HardyRama}).  
		Then $\Phi_c$ extends $\Omega$ if and only if $c(\rat,p) = -1/(\log p)$ for all $p\ne \infty$.
	\item[Integer Logarithm] Let $\Psi(n)$ be the sum of the prime factors of $n$, counted with multiplicity (see \cite{Jaki,Lal,Atanassov}).  Then $\Phi_c$ extends $\Psi$ if and only if 
		$c(\rat,p) = -p/(\log p)$ for all $p\ne \infty$.
\end{description}
It would be interesting to discover a version of Corollary \ref{RestrictedIso} that could be applied to all additive functions rather than only to completely additive functions.  For now, we are unaware of
any way to formulate such a result.

Of the previous work on this subject, we should regard Theorem \ref{RealSurjective} as most analogous to \cite[Theorem 1.3]{SamuelsConsistent}, which established the right hand isomorphism of \eqref{LCIsos}.
It is important to note, however, that our result cannot be proved by directly applying existing work.  Although $\mathcal G$ appears as a dense subset of $LC_c^0(Y)$ with respect to the 
$L^1$-norm, we impose no continuity assumption on elements of $\mathcal L(\mathcal G,\real)$.  Hence, prior to proving Theorem \ref{RealSurjective},
we cannot be certain that an arbitrary map $\Phi\in \mathcal L(\mathcal G,\real)$ may be extended to a linear map on $LC_c^0(Y)$.  While our proof is inspired by some ideas of 
\cite{SamuelsConsistent}, several different methods are needed, including an application of Dirichlet's $S$-unit Theorem.

Related to these observations, there is a version of Theorem \ref{RealSurjective} which classifies all continuous linear maps from $\mathcal G$ to $\real$, where $\mathcal G$ is equipped with the Weil height norm as in
\cite{AllcockVaaler}.  Specifically, if we let $\mathcal J'$ be the set of consistent maps $c$ for which
\begin{equation*}
	\sup\left\{\left|\frac{c(K,v)}{\lambda(K,v)}\right|: (K,v)\in \mathcal J\right\} < \infty,
\end{equation*}
then $c\mapsto \Phi_c$ is an isomorphism of $\mathcal J'$ onto the space of continuous linear maps from $\mathcal G$ to $\real$.
However, this fact does indeed follow directly from previous work, namely \cite[Theorem 1.6]{SamuelsConsistent}.  Because of our observations following Corollary \ref{RestrictedIso}, both the natural logarithm and the 
prime Omega function have continuous extensions to $\mathcal G$, while the integer logarithm has no such extension.

\subsection{Rational Valued Linear Maps} With an eye toward studying $\mathcal G^*$, we are particularly interested in identifying consistent maps $c$ for which $\Phi_c(\alpha) \in \rat$ for all $\alpha\in \mathcal G$.  
To facilitate these efforts, we shall write
\begin{equation*}
	\mathcal I^* = \left\{c\in \mathcal J^*: \Phi^*(c) \in \mathcal G^*\right\} = \left\{c\in \mathcal J^*: \Phi_c(\alpha) \in \rat\mbox{ for all }\alpha\in \mathcal G\right\} 
\end{equation*}
and note that $\mathcal I^*$ is a $\rat$-linear subspace of $\mathcal J^*$.  The question now arises to provide a necessary and sufficient condition for the claim that $c\in \mathcal I^*$.  This problem appears to be quite 
challenging, however, we may obtain several interesting examples by applying two supplementary results, the first of which is as follows. 

\begin{thm} \label{CanonicalConsistent}
	Suppose $K$ is a number field, and for each place $v$ of $K$, let $y_v\in \real$.  Then there exists a unique consistent map $c\in \mathcal J^*$ such that
	\begin{equation} \label{DF}
		c(L,w) = \frac{[L_w:K_v]}{[L:K]} y_v
	\end{equation}
	for all $v\in M_K$, all finite extensions $L/K$, and all places $w$ of $L$ dividing $v$.  Moreover, if $\Phi_c(\alpha) \in \rat$ for all $\alpha\in K^\times$ then $c\in \mathcal I^*$.	
\end{thm}

Given a point $\yy = (y_v)_{v\in M_K}$, we write $c_\yy$ to denote the consistent map from Theorem \ref{CanonicalConsistent}.  We plainly have that $c_\yy(K,v) = y_v$ for all places $v$ of $K$,
and therefore, if we wish for $c_\yy\in \mathcal I^*$, it is sufficient to find a point $\yy = (y_v)_{v\in M_K}$ such that
\begin{equation*}
	\sum_{v\in M_K} y_v\log\|\alpha\|_v\in \rat \quad\mbox{for all }\alpha\in K^\times.
\end{equation*}
In the special case where $K =\rat$ or where $K$ an imaginary quadratic extension of $\rat$, we can locate examples of this sort rather easily.  Indeed, such number fields have a unique 
Archimedean place $\infty$, so we may choose $y_\infty = 0$ and $y_v\log p_v \in \rat$ for all non-Archimedean places $v$ of $K$.  The resulting consistent map $c_\yy$ certainly belongs to $\mathcal I^*$.

The more interesting cases, however, arise from number fields having multiple Archimedean places.  The following theorem is somewhat technical, but useful in constructing other examples of consistent maps in $\mathcal I^*$.

\begin{thm}\label{KRational}
	Suppose $K$ is a number field.  Let $\{v_1,v_2,\ldots,v_n\}$ be the complete list of Archimedean places of $K$ and let $\{\alpha_1,\alpha_2,\ldots,\alpha_{n-1}\}$ be a fundamental set of units in $K$.  
	For each non-Archimedean place $v$ of $K$, let $\beta_v\in \mathcal O_K$ be such that $\|\beta_v\|_v < 1$ and $\|\beta_v\|_w = 1$ for all non-Archimedean places $w\ne v$.  Finally, let $\yy = (y_v)_{v\in M_K}$.  Then 
	\begin{equation*}
		\sum_{v\in M_K} y_v\log\|\alpha\|_v\in \rat \quad\mbox{for all }\alpha\in K^\times
	\end{equation*}
	if and only if the following two conditions hold:
	\begin{enumerate}[label=(\roman*)]
		\item \label{Arches}We have
			\begin{equation*}
				\begin{pmatrix}
					\log\|\alpha_1\|_{v_1} & \log\|\alpha_1\|_{v_2} & \cdots & \log\|\alpha_1\|_{v_n} \\
					\log\|\alpha_2\|_{v_1} & \log\|\alpha_2\|_{v_2} & \cdots & \log\|\alpha_2\|_{v_n} \\
					\vdots & \vdots & \ddots & \vdots \\
					\log\|\alpha_{n-1}\|_{v_1} & \log\|\alpha_{n-1}\|_{v_2} & \cdots & \log\|\alpha_{n-1}\|_{v_n}
				\end{pmatrix}
				\begin{pmatrix} y_{v_1} \\ y_{v_2} \\ \vdots \\ y_{v_n}\end{pmatrix} \in \rat^{n-1}
			\end{equation*}
		\item\label{NonArches} For all non-Archimedean places $v$ of $K$, we have that
			\begin{equation*}
				y_v\log\|\beta_v\|_v + \sum_{i=1}^n y_{v_i}\log \|\beta_v\|_{v_i} \in \rat.
			\end{equation*}
	\end{enumerate}
\end{thm}

The existence of the fundamental set of units $\{\alpha_1,\alpha_2,\ldots,\alpha_{n-1}\}$ is guaranteed by Dirichlet's Unit Theorem, and furthermore, \cite[Lemma 3.1]{SamuelsClassification} establishes the existence of the
points $\beta_v$.  As a result, it is always possible to select points $\alpha_i$ and $\beta_v$ that satisfy the assumptions of Theorem \ref{KRational}.  Now we may use Theorems \ref{CanonicalConsistent}
and \ref{KRational} to construct maps $c\in \mathcal I^*$ which are distinct from known examples.

\begin{ex}
	We consider $K = \rat(\sqrt 2)$ and fix an embedding of $K$ into $\real$ with $|\cdot |$ denoting the usual absolute value on $\real$.  We know that $K$ has two Archimedean places which we shall denote by
	$v_1$ and $v_2$.  We may assume that
	\begin{equation*}
		\|a+b\sqrt 2\|_{v_1} = |a + b\sqrt 2|\quad\mbox{and}\quad \|a+b\sqrt 2\|_{v_2} = |a-b\sqrt 2|
	\end{equation*}
	for all $a,b\in \intg$.  We observe that $K$ has class number equal to $1$, its ring of integers is $\mathcal O_K = \intg[\sqrt 2]$,
	and $1+\sqrt 2$ is the fundamental unit in $K$.  Additionally, a prime $p\in \rat$ splits in $K$ if and only if $p\equiv 1, 7 \mod 8$.
	
	We now select the points $y_v$ for use in Theorem \ref{KRational}.
	\begin{enumerate}[label=(\roman*)]
		\item We define
			\begin{equation*}
				y_{v_1} = \frac{1}{\log (1+\sqrt 2)}  \quad\mbox{and}\quad y_{v_2} = -\frac{1}{\log (1+\sqrt 2)}.
			\end{equation*}
		\item If $p \not\equiv 1, 7\mod 8$ and $v$ divides $p$, we let $y_v = 0$.
		\item\label{NonArchPlaces} If $p\equiv 1, 7\mod 8$ then we let $v$ and $w$ be distinct places of $K$ dividing $p$.  In this case, $p$ has the form $p=\beta_v\beta_w$, where $\beta_v$ and $\beta_w$ are 
			generators of the prime ideals
			\begin{equation*}
				\left\{\alpha\in\mathcal O_K: \|\alpha\|_v < 1\right\}\quad\mbox{and}\quad \left\{\alpha\in\mathcal O_K: \|\alpha\|_w < 1\right\},
			\end{equation*}
			respectively.  In this situation, we define
			\begin{equation*}
				y_v = \frac{\log\|\beta_v\|_{v_1} - \log\|\beta_v\|_{v_2}}{(\log p)(\log(1+\sqrt 2))}\quad\mbox{\and}\quad y_w = \frac{\log\|\beta_w\|_{v_1} - \log\|\beta_w\|_{v_2}}{(\log p)(\log(1+\sqrt 2))}.
			\end{equation*}
	\end{enumerate}
	It is straightforward to verify that the points $y_v$ satisfy the two properties of Theorem \ref{KRational}, and then by Theorem \ref{CanonicalConsistent}, $c = c_\yy\in \mathcal I^*$.   In other words, $\Phi_c$
	is a (rational-valued) linear functional on $\mathcal G$.
	
	Since $\beta_v$ and $\beta_w$ are conjugates over $\rat$, the values in \ref{NonArchPlaces} satsify $y_v = -y_w$.  Therefore, we have $c(K,v) = -c(K,w)$ whenever $v$ and $w$ divide the same place of $\rat$, 
	or equivalently, we have $c(\rat,p)=0$ for all $p\in M_\rat$.  The first few non-zero values of $c(K,v)$ are approximated in the following table.
	
	\begin{equation*}
	\begin{array}{|c|c|c|}
		p & \mbox{Factorization of } p\mbox{ in } \intg[\sqrt 2] & c(K,v)\mbox{ for } v\mid p\mbox{ (approx.)}\\ \hline \hline
		\infty & \mbox{NA} &  \pm 1.13459 \\ \hline
		7 & (3+\sqrt{2})(3-\sqrt 2) & \pm 0.596913 \\ \hline
		17 & (5 + 2\sqrt 2)(5-2\sqrt 2)& \pm 0.513516  \\ \hline
		23 & (5 + \sqrt 2)(5 - \sqrt 2) & \pm 0.0.513516 \\ \hline
		31 & (7 + 3\sqrt 2)(7 - 3\sqrt 2) & \pm 0.464359 \\ \hline
		41 & (7 + 2\sqrt 2)(7 - 2\sqrt 2) &\pm 0.261831 \\ \hline
		47 & (7 + \sqrt 2)(7-\sqrt 2) & \pm 0.120733 \\ \hline
		71 & (11 + 5\sqrt 2)(11-5\sqrt 2) & \pm 0.406159 \\ \hline
	\end{array}
	\end{equation*}
\end{ex}
\smallskip

\subsection{Organizational Summary} We shall structure the remainder of this article by separating the proof of Theorem \ref{RealSurjective} into two components.  First, in Section \ref{Inj}, we show that $\Phi^*$ is a linear transformation 
such that $\ker(\Phi^*) = \spann_\real\{\lambda\}$.   The surjectivity component of the proof requires applying Theorem \ref{CanonicalConsistent}, and hence, we use Section \ref{Canonical} to prove that result.
Finally, in Section \ref{Surj}, we complete the proof of Theorem \ref{RealSurjective} by proving that $\Phi^*$ is surjective.  The proof of Theorem \ref{KRational} is included in that section as well. 

\section{The Kernel of $\Phi^*$} \label{Inj}

\begin{thm} \label{KernelOnly}
	The map $\Phi^*:\mathcal J^* \to \mathcal L(\mathcal G,\real)$ is an $\real$-linear transformation such that $\ker(\Phi^*) = \spann_\real\{\lambda\}$.
\end{thm}

Before proceeding with the proof of Theorem \ref{KernelOnly}, we remind the reader of the relevant features of Dirichlet's Unit Theorem \cite[Theorem 7.31]{Jarvis}.
If $K$ is a number field, then $\mathcal O_K$ denotes its ring of integers and 
\begin{equation*}
	\mathcal O_K^\times = \left\{\alpha\in K: \|\alpha\|_v = 1\mbox{ for all } v\nmid\infty\right\}
\end{equation*}
is called its {\it group of units}.  If $K$ has $n$ Archimedean places, then Dirchlet's Unit Theorem asserts that the $\mathcal O_K^\times$ has rank equal to $n-1$.
If $\zeta$ is a root of unity and $\alpha_1,\alpha_2,\ldots,\alpha_{n-1}\in \mathcal O_K$ are such that $\{\zeta,\alpha_1,\alpha_2,\cdots,\alpha_{n-1}\}$ generates $\mathcal O_K^\times$, 
then the collection $\{\alpha_1,\alpha_2,\cdots,\alpha_{n-1}\}$ is called a {\it fundamental set of units} in $K$.

The proof of Theorem \ref{KernelOnly} begins with the following lemma.

\begin{lem} \label{RankRegulator}
	Let $K$ be a number field having Archimedean places $\{v_1,v_2,\ldots,v_n\}$, let $\{\alpha_1,\alpha_2,\ldots,\alpha_{n-1}\}$ be a fundamental set of units in $K$, and define
	\begin{equation*}
		A = \begin{pmatrix}
				\log\|\alpha_1\|_{v_1} & \log\|\alpha_1\|_{v_2} & \cdots & \log\|\alpha_1\|_{v_n} \\
				\log\|\alpha_2\|_{v_1} & \log\|\alpha_2\|_{v_2} & \cdots & \log\|\alpha_2\|_{v_n} \\
				\vdots & \vdots & \ddots & \vdots \\
				\log\|\alpha_{n-1}\|_{v_1} & \log\|\alpha_{n-1}\|_{v_2} & \cdots & \log\|\alpha_{n-1}\|_{v_n}
			\end{pmatrix}.
	\end{equation*}
	Then $\rank(A) = n-1$ and $\dim(\ker A) = 1$.
\end{lem}
\begin{proof}
	For simplicity, we write $D_i = [K_{v_i}:\rat_{v_i}]$ and define the following additional matrices:
	\begin{equation*}
		D = \begin{pmatrix} D_1 & 0 & \cdots & 0 \\
						0 & D_2 & \cdots & 0 \\
						\vdots & \vdots & \ddots & \vdots \\
						0 & 0 & \cdots & D_n \end{pmatrix}
	\end{equation*}
	and 
	\begin{align*}
		B & = AD = \begin{pmatrix}
				D_1\log\|\alpha_1\|_{v_1} & D_2\log\|\alpha_1\|_{v_2} & \cdots & D_n\log\|\alpha_1\|_{v_n} \\
				D_1\log\|\alpha_2\|_{v_1} & D_2\log\|\alpha_2\|_{v_2} & \cdots & D_n\log\|\alpha_2\|_{v_n} \\
				\vdots & \vdots & \ddots & \vdots \\
				D_1\log\|\alpha_{n-1}\|_{v_1} & D_2\log\|\alpha_{n-1}\|_{v_2} & \cdots & D_n\log\|\alpha_{n-1}\|_{v_n}
			\end{pmatrix}.
	\end{align*}
	Clearly $\det(D)\ne 0$ so that $A$ and $B$ must have the same rank.  However, if we let $B_i$ denote the matrix obtained by removing column $i$ from $B$, then it is well-known that
	$|\det(B_i)|$ is non-zero and independent of $i$.  This value is called the {\it regulator of $K$} and is thoroughly studied throughout the literature on algebraic number theory (see \cite[Def. 10.8]{Jarvis}, for example).
	In any case, it now follows that the rows of $B$ are linearly independent so that
	\begin{equation*}
		\rank(A) = \rank(B) = n-1.
	\end{equation*}
	Now applying the rank-nullity theorem, we conclude that $\dim(\ker(A)) = 1$.
\end{proof}

Our next result is the primary ingredient in identifying the kernel of $\Phi^*$.   As we shall find, it also plays a crucial role in showing that $\Phi^*$ is surjective.

\begin{lem}\label{ZeroPhi}
	Suppose that $c:\mathcal J\to \real$ is a consistent map and $K$ is a number field.  If $\Phi_c(\alpha) = 0$ for all $\alpha\in K^\times$ then 
	\begin{equation} \label{InSpanLambda}
		c(K,v) = c(\rat,\infty) \lambda(K,v)
	\end{equation}
	for all $v\in M_K$.
\end{lem}
\begin{proof}
	We first establish \eqref{InSpanLambda} in the case that $v$ is Archimedean.  To this end, we let $v_1,v_2,\ldots,v_n$ be the complete list of Archimedean places of $K$ so that the consistency property
	\eqref{Consistency} implies that
	\begin{equation} \label{ArchCons}
		c(\rat,\infty) = \sum_{i=1}^n c(K,v_i).
	\end{equation}
	If $n=1$ then $\lambda(K,v_1) = 1$, and the required property follows immediately from \eqref{ArchCons}.  Therefore, we shall assume that $n\geq 2$.
	
	According to Dirichlet's Unit Theorem, we may let $\{\alpha_1,\alpha_2,\ldots,\alpha_{n-1}\}$ be a
	set of fundamental units in $K$ and let $A$ be the $(n-1)\times n$ matrix $A$ given by Lemma \ref{RankRegulator}.  Additionally, we let
	\begin{equation*}
		C = \begin{pmatrix} c(K,v_1) \\ c(K,v_2) \\ \vdots \\ c(K,v_n) \end{pmatrix}\quad\mbox{and}\quad \Lambda =  \begin{pmatrix} \lambda(K,v_1) \\ \lambda(K,v_2) \\ \vdots \\ \lambda(K,v_n) \end{pmatrix}.
	\end{equation*}
	As we have assumed that $\Phi_c(\alpha) = 0$ for all $\alpha\in K$, we have
	\begin{equation} \label{CKernel}
		AC = \begin{pmatrix} \Phi_c(\alpha_1) \\ \Phi_c(\alpha_2) \\ \vdots \\ \Phi_c(\alpha_{n-1}) \end{pmatrix} = \begin{pmatrix} 0 \\ 0 \\ \vdots \\ 0 \end{pmatrix}
	\end{equation}
	so that $C\in \ker(A)$.  However, the product formula implies that $\Lambda$ also belongs to $\ker(A)$, and since $\Lambda$ is clearly non-zero, Lemma \ref{RankRegulator} establishes that
	$\ker(A) = \left\{ r\Lambda: r\in \real\right\}$.   We now obtain a real number $r$ such that
	\begin{equation*}
		c(K,v_i) = r \lambda(K,v_i)\quad\mbox{ for all } 1\leq i\leq n.
	\end{equation*}
	Finally, property \eqref{ArchCons} shows that
	\begin{equation*}
		c(\rat,\infty) = \sum_{i=1}^n c(K,v_i) = r \sum_{i=1}^n \lambda(K,v_i) = r
	\end{equation*}
	and we have established \eqref{InSpanLambda} for all Archimedean places $v$ of $K$.

	We now establish \eqref{InSpanLambda} when $v$ is non-Archimedean.  Because of what we have already shown, we may assume for the remainder of this proof that 
	\begin{equation} \label{QBaseArch}
		c(K,u) =  c(\rat,\infty) \frac{[K_u:\rat_\infty]}{[K:\rat]}
	\end{equation}
	for all Archimedean places $u$ of $K$.  According to \cite[Theorem 3.1]{SamuelsClassification}, there exists $\beta\in K$ such that $\|\beta\|_v < 1$ and $\|\beta\|_w = 1$ for all other non-Archimedean places $w$ of $K$.
	Now applying \eqref{QBaseArch}, we obtain
	\begin{align*}
		\Phi_c(\beta) & = c(K,v)\log \|\beta\|_v + \sum_{u\mid\infty} c(K,u) \log\|\beta\|_u \\
		& = c(K,v)\log \|\beta\|_v + c(\rat,\infty) \sum_{u\mid\infty} \frac{[K_u:\rat_\infty]}{[K:\rat]} \log\|\beta\|_u.
	\end{align*}
	According to the product formula on $K$, the summation on the right hand side may be simplified so that
	\begin{equation*}
		\Phi_c(\beta) = c(K,v)\log \|\beta\|_v - c(\rat,\infty) \frac{[K_v:\rat_v]}{[K:\rat]} \log\|\beta\|_v.
	\end{equation*}
	From our assumptions we have that $\Phi_c(\beta) = 0$ and $\|\beta\|_v \ne 1$, so it follows that
	\begin{equation*}
		c(K,v) = c(\rat,\infty) \frac{[K_v:\rat_v]}{[K:\rat]}
	\end{equation*}
	as required.
\end{proof}

\begin{proof}[Proof of Theorem \ref{KernelOnly}]
	Suppose that $c,d\in \mathcal J^*$ and $r\in \real$.  For each $\alpha\in \algt$, we assume that $K$ is a number field containing $\alpha$ and observe that
	\begin{align*}
		\Phi_{c+d}(\alpha) & = \sum_{v\in M_K} [c(K,v) + d(K,v)]\log\|\alpha\|_v \\
					& = \sum_{v\in M_K} c(K,v)\log\|\alpha\|_v + \sum_{v\in M_K} d(K,v)\log\|\alpha\|_v \\
					& = \Phi_c(\alpha) + \Phi_d(\alpha)
	\end{align*}
	which proves that $\Phi^*(c+d) = \Phi^*(c) + \Phi^*(d)$.  Also, we have
	\begin{equation*}
		\Phi_{rc}(\alpha) = r\sum_{v\in M_K} c(K,v)\log\|\alpha\|_v = r\Phi_c(\alpha)
	\end{equation*}
	establishing that $\Phi^*(rc) = r\Phi^*(c)$ and showing that $\Phi^*$ is a linear transformation.

	Assuming that $c\in \spann_\real\{\lambda\}$, there exists $r\in \real$ such that 
	\begin{equation*}
		c(K,v) = r\cdot \frac{[K_v:\rat_v]}{[K:\rat]}
	\end{equation*}
	for all number fields $K$ and all places $v$ of $K$.  For each non-zero point $\alpha\in \alg$, the product formula now implies that
	\begin{align*}
		\Phi_c(\alpha) & = \sum_{v\in M_K} c(K,v)\log\|\alpha\|_v  = r \sum_{v\in M_K} \frac{[K_v:\rat_v]}{[K:\rat]}\log \|\alpha\|_v = 0 
	\end{align*}
	proving that $\Phi_c\equiv 0$ and $\spann_\real\{\lambda\}\subseteq \ker(\Phi^*)$.
	
	Now assuming that $\Phi_c\equiv 0$ and $(K,v)\in \mathcal J$, we certainly have that $\Phi_c(\alpha) = 0$ for all $\alpha\in K^\times$.  Hence, Lemma \ref{ZeroPhi} applies to yield
	\begin{equation*}
		c(K,v) = c(\rat,\infty)\lambda(K,v)
	\end{equation*}
	for all $v\in M_K$ establishing that $c\in \spann_\real\{\lambda\}$ and $\ker(\Phi^*) \subseteq \spann_\real\{\lambda\}$, as required.
\end{proof}

\section{Extensions of Consistent Maps} \label{Canonical}

Before we continue with the proof of Theorem \ref{RealSurjective}, we will provide our proof of Theorem \ref{CanonicalConsistent}.  As we shall find, Theorem \ref{CanonicalConsistent} is required in the proof of Theorem 
\ref{RealSurjective}, and as such, it makes sense to provide its proof first.  

If $K$ is a number field and $L/K$ is a finite extension, we note the two well-known identities
\begin{equation} \label{LocalGlobal}
	[L:K] = \sum_{w\mid v} [L_w:K_v]\quad\mbox{and}\quad \norm_{L/K}(\alpha) = \prod_{w\mid v} \norm_{L_w/K_v}(\alpha)
\end{equation}
for all $\alpha\in L$ (see \cite[Eq. (2) and Prop. 4]{Frohlich}, for example).

\begin{proof}[Proof of Theorem \ref{CanonicalConsistent}]
	If $L$ is any number field, let $F$ be a finite extension of $L$ containing $K$.  Further, if $t$ is a place of $F$ dividing the place $v$ of $K$, then we write $y_t = y_v$.  For each $w\in M_L$ we define
	\begin{equation} \label{DFDef}
		d_F(L,w) = \sum_{t\mid w} \frac{[F_t:K_t]}{[F:K]} y_t,
	\end{equation}
	where the summation on the right hand side of \eqref{DFDef} runs over places $t$ of $F$ dividing $w$.  We claim that $d_F(L,w)$ is independent of $F$.  
	
	To see this, we suppose that $E$ is a finite extension of $F$.  In the ensuing calculations, we shall adopt the convention of writing $t$ for places of $F$ and $s$ for places of $E$.
	Clearly if $s\mid t$ then $y_t = y_s$ and $K_s = K_t$, so we obtain
	\begin{align*}
		d_E(L,w) & = \sum_{s\mid w} \frac{[E_s:K_s]}{[E:K]} y_s \\
				& = \sum_{t\mid w} \sum_{s\mid t} \frac{[E_s:F_t]\cdot[F_t:K_t]}{[E:F]\cdot[F:K]} y_t \\
				& = \sum_{t\mid w} \frac{[F_t:K_t]}{[F:K]} y_t  \sum_{s\mid t} \frac{[E_s:F_t]}{[E:F]}.
	\end{align*}
	Now applying the left hand equality of \eqref{LocalGlobal}, we obtain that
	\begin{equation*}
		d_E(L,w) = \sum_{t\mid w} \frac{[F_t:K_t]}{[F:K]} y_t = d_F(L,w)
	\end{equation*}
	showing that $d_F(L,w)$ is indeed independent of $F$.  Hence, we may define $c:\mathcal J\to \real$ by $c(L,w) = d_F(L,w)$, where $F$ is any number field containing both $K$ and $L$.
	
	To prove that $c$ is consistent, we assume that $M$ is a finite extension of $L$ and $w$ is a place of $L$.  We select a number field $F$ containing both $M$ and $K$.
	Then using $x$ to denote places of $M$, we obtain from \eqref{DFDef} that
	\begin{align*}
		\sum_{x\mid w} c(M,x) & = \sum_{x\mid w} d_F(M,x) \\
						& = \sum_{x\mid w} \sum_{t\mid x} \frac{[F_t:K_t]}{[F:K]} y_t \\
						& = \sum_{t\mid w} \frac{[F_t:K_t]}{[F:K]} y_t  = d_F(L,w)  = c(L,w),
	\end{align*}
	proving that $c$ is consistent.  
	
	To establish \eqref{DF}, we assume that $L$ is a finite extension of $K$, $v$ is a place of $K$, and $w$ is a a place of $L$ dividing $v$.  Now apply \eqref{DFDef} with $F=L$
	to obtain
	\begin{equation*}
		c(L,w) =\frac{[L_w:K_w]}{[L:K]}y_w = \frac{[L_w:K_v]}{[L:K]}y_v
	\end{equation*}
	which is the required property \eqref{DF}.
	
	To prove that this consistent map is unique, we suppose that $c,d\in \mathcal J^*$ both satisfy \eqref{DF}.  This means $c(L,w) = d(L,w)$ for all finite extensions $L/K$ and all places $w$ of $L$.
	Now if $L'$ is an arbitrary number field and $w'$ is a place of $L'$, we may choose a number field $L$ containing both $K$ and $L'$.  Then applying the consistency of $c$ and $d$, we obtain
	\begin{equation*}
		c(L',w') = \sum_{w\mid w'} c(L,w) = \sum_{w\mid w'} d(L,w) = d(L',w')
	\end{equation*}
	proving that $c= d$.
	
	Finally, we assume that $\Phi_c(\alpha) \in \rat$ for all $\alpha\in K^\times$.  If $\beta\in \algt$ we may let $L$ be a number field containing both $\beta$ and $K$.  Now applying \eqref{DF} we obtain
	\begin{align*}
		\Phi_c(\beta) & = \sum_{w\in M_L} c(L,w) \log \|\beta\|_w \\
					& = \sum_{v\in M_K} \sum_{w\mid v} c(L,w) \log \|\beta\|_w = \sum_{v\in M_K} y_v \sum_{w\mid v} \frac{[L_w:K_v]}{[L:K]}\log \|\beta\|_w.
	\end{align*}
	Then using the right hand equality of \eqref{LocalGlobal}, we find that
	\begin{align*} 
			\Phi_c(\beta) & = \sum_{v\in M_K} y_v \sum_{w\mid v} \frac{[L_w:K_v]}{[L:K]}\log \|\norm_{L_w/K_v}(\beta)\|_v^{1/[L_w:K_v]} \\
					& = \frac{1}{[L:K]}  \sum_{v\in M_K} y_v \sum_{w\mid v}\log \|\norm_{L_w/K_v}(\beta)\|_v \\
					& = \frac{1}{[L:K]}  \sum_{v\in M_K} y_v\log \|\norm_{L/K}(\beta)\|_v \\
					& =  \frac{1}{[L:K]} \Phi_c(\norm_{L/K}(\beta)).
	\end{align*}
	which is clearly rational.
\end{proof}

\section{Surjectivity of $\Phi^*$} \label{Surj}

In order to complete the proof of Theorem \ref{RealSurjective}, we must establish the following result.

\begin{thm} \label{SurjectiveOnly}
	The map $\Phi^*:\mathcal J^*\to \mathcal L(\mathcal G,\real)$ is surjective.
\end{thm}

The proof of Theorem \ref{SurjectiveOnly} requires some background notation as well as three lemmas.  Let $K$ be a number field and let $S$ be a finite subset of $M_K$ containing the Archimedean 
places of $K$.  The set
\begin{equation*}
	U_{K,S} = \left\{ \alpha\in K^\times: \|\alpha\|_v = 1\mbox{ for all } v\in M_K\setminus S\right\}
\end{equation*}
is a subgroup of $K^\times$ called the {\it group of $S$-units in $K$}, which according to Dirichlet's $S$-unit Theorem (see \cite[Thm. III.3.5]{Nark} or \cite[\S 1.1]{FiliMinerDirichlet}),
is finitely generated of rank $\#S - 1$.   If $\zeta$ is a root of unity and $\alpha_1,\alpha_2,\ldots,\alpha_{\#S-1}\in U_{K,S}$ are such that $\{\zeta,\alpha_1,\alpha_2,\cdots,\alpha_{\#S-1}\}$ generates $U_{K,S}$, 
the then collection $\{\alpha_1,\alpha_2,\cdots,\alpha_{\#S-1}\}$ is called a {\it fundamental set of $S$-units} in $K$.  Of course, this terminology generalizes our definitions from the beginning of Section \ref{Inj}.

For our proof that $\Phi^*$ is surjective, we require a reinterpretation of Dirichlet's $S$-unit Theorem in the language of linear algebra.  To this end, we write
\begin{equation*}
	K_\divv = \left\{\alpha\in \algt:\alpha^n\in K\mbox{ for some } n\in \nat\right\},
\end{equation*}
and note that $K_\divv$ is a subgroup of $\algt$ containing $\algt_\tors$. We further write
\begin{equation*}
	\mathcal G_K = K_\divv/\algt_\tors
\end{equation*}
and observe that $\mathcal G_K$ is a subspace of $\mathcal G$.  If $\alpha\in K_\divv$, we shall adopt the convention of writing $\bar\alpha$ to denote its image under the canonical homomorphism 
$K_\divv\to K_\divv/\algt_\tors$.  For each place $v$ of $K$ and $\alpha\in K_\divv$, we may assume that $\ell\in \nat$ is such that $\alpha^\ell\in K$ and define
\begin{equation*}
	\|\alpha\|_v = \|\alpha^\ell\|_v^{1/\ell}.
\end{equation*}
The right hand side of this equality does not depend on the choice of $\ell$, and hence, $\|\cdot\|_v$ is a well-defined map on $\mathcal G_K$ which satisfies
\begin{enumerate}[label=(A\arabic*)]
	\item\label{Multiply} $\|\alpha\beta\|_v = \|\alpha\|_v\|\beta\|_v$ for all $\alpha,\beta\in K_\divv$
	\item\label{Scalar} $\|\alpha^r\|_v = \|\alpha\|_v^r$ for all $\alpha\in K_\divv$ and all $r\in \rat$.
\end{enumerate}
Equivalently, $\alpha\mapsto \log\|\alpha\|_v$ defines a linear transformation from $\mathcal G_K$ to $\real$ when viewed as $\rat$-vector spaces.  Finally, we let
\begin{equation*}
	\mathcal G_{K,S} = \left\{ \alpha\in \mathcal G_K: \|\alpha\|_v = 1\mbox{ for all } v\in M_K\setminus S\right\}
\end{equation*}
and note the following manner of identifying a basis for $\mathcal G_{K,S}$ over $\rat$.

\begin{lem} \label{Dirichlet}
	Suppose that $K$ is a number field having $n$ Archimedean places.  Assume that $S_\infty$ is the complete set of Archimedean places of $K$, $S_0 = \{w_1,w_2,\ldots,w_m\}$ is a finite (possibly empty) set of 
	non-Archimedean places of $K$, and $S = S_\infty \cup S_0$.  Further assume the following:
	\begin{enumerate}[label=(\Roman*)]
		\item\label{FundUnits} $\{\alpha_1,\alpha_2,\ldots,\alpha_{n-1}\}$ is a fundamental set of units in $K$
		\item\label{NonArchChoices} For every $1\leq i\leq m$, $\beta_i\in K$ is such that $\|\beta_i\|_{w_i} < 1$ and $\|\beta_i\|_{w} = 1$ for all non-Archimedean places $w\ne w_i$.
	\end{enumerate}
	Then $\{\bar\alpha_1,\bar\alpha_2,\ldots,\bar\alpha_{n-1},\bar\beta_1,\bar\beta_2,\ldots,\bar\beta_m\}$ is a basis for $\mathcal G_{K,S}$ over $\rat$.  In particular,
	$\mathcal G_{K,S}$ is a finite dimensional subspace of $\mathcal G_K$ with $\dim(\mathcal G_{K,S}) = \#S-1$.
\end{lem}
\begin{proof}
	For each $0\leq k\leq m$, we define $T_k = S_\infty\cup\{w_1,w_2,\ldots,w_k\}$ so that $T_m = S$.  We shall prove by induction on $k$ that
	\begin{equation*}
		\{\bar\alpha_1,\bar\alpha_2,\ldots,\bar\alpha_{n-1},\bar\beta_1,\bar\beta_2,\ldots,\bar\beta_k\}
	\end{equation*}
	is a basis for $\mathcal G_{K,T_k}$ for all $0\leq k \leq m$.  The lemma would then follow by taking the special case $k = m$.
		
	{\bf Base Case:} Since we have assumed that $\{\alpha_1,\alpha_2,\ldots,\alpha_{n-1}\}$ is a fundamental set of units in $K$, we obtain that
	\begin{enumerate}[label=(\roman*)]
		\item\label{Spans} For every $\alpha\in U_{K,T_0}$, there exist $r_1,r_2,\ldots,r_{n-1}\in \intg$ such that $\alpha = \alpha_1^{r_1}\alpha_2^{r_2}\cdots \alpha_{n-1}^{r_{n-1}}$.
		\item\label{LI} If $r_1,r_2,\ldots,r_{n-1}\in \intg$ are such that $\alpha_1^{r_1}\alpha_2^{r_2}\cdots \alpha_{n-1}^{r_{n-1}} = 1$ then $r_i=0$ for all $i$.
	\end{enumerate}
	Let $\bar\alpha\in \mathcal G_{K,T_0}$ and assume that $\alpha$ is some representative of $\bar\alpha$ in $K_\divv$.  By definition of $K_\divv$, there exists $\ell\in \nat$ such that $\alpha^\ell\in K^\times$.  
	But clearly we also have that $\|\alpha^\ell\|_w = 1$ for all non-Archimedean places $w$ of $K$, and hence, $\alpha^\ell\in U_{K,T_0} = U_{K,S_\infty}$.  Then according to \ref{Spans}, there exist $r_1,r_2,\ldots,r_{n-1}\in \intg$ 
	such that $\alpha^\ell = \alpha_1^{r_1}\alpha_2^{r_2}\cdots \alpha_{n-1}^{r_{n-1}}$.  We now conclude that
	\begin{equation*}
		\bar\alpha^\ell = \bar\alpha_1^{r_1}\bar\alpha_2^{r_2}\cdots \bar\alpha_{n-1}^{r_{n-1}}
	\end{equation*}
	Since $\ell$ is a positive integer, we get that
	\begin{equation*}
		\bar\alpha= \bar\alpha_1^{r_1/\ell}\bar\alpha_2^{r_2/\ell}\cdots \bar\alpha_{n-1}^{r_{n-1}/\ell},
	\end{equation*}
	proving that $\bar\alpha\in \spann\{\bar\alpha_1,\bar\alpha_2,\ldots,\bar\alpha_{n-1}\}$.
	
	Now assume that $a_i\in \intg$ and $b_i\in \nat$ are such that
	\begin{equation*}
		\bar\alpha_1^{a_1/b_1}\bar\alpha_2^{a_2/b_2}\cdots \bar\alpha_{n-1}^{a_{n-1}/b_{n-1}} = 1.
	\end{equation*}
	After raising both sides to the $b_1b_2\cdots b_{n-1}$ power, we obtain
	\begin{equation*}
		1 = \bar\alpha_1^{r_1}\bar\alpha_2^{r_2}\cdots \bar\alpha_{n-1}^{r_{n-1}} = \overline{ \alpha_1^{r_1}\alpha_2^{r_2}\cdots \alpha_{n-1}^{r_{n-1}} },\quad \mbox{where } r_i = a_i\prod_{j\ne i} b_i
	\end{equation*}
	As a result, there must exist a root of unity $\zeta$ such that $\zeta =  \alpha_1^{r_1}\alpha_2^{r_2}\cdots \alpha_{n-1}^{r_{n-1}}.$
	We certainly have that $\zeta^d = 1$ for some $d\in \nat$, and hence, 
	\begin{equation*}
		1 =  \alpha_1^{dr_1}\alpha_2^{dr_2}\cdots \alpha_{n-1}^{dr_{n-1}}.
	\end{equation*}
	By applying \ref{LI}, we conclude that $dr_i=0$ for all $i$, and since $d$ is certainly non-zero, we obtain that $r_i=0$.  It now follows that $a_i=0$, as required.
	
	{\bf Inductive Step}: We now let
	\begin{equation*}
		\mathcal B_k = \{\bar\alpha_1,\bar\alpha_2,\ldots,\bar\alpha_{n-1},\bar\beta_1,\bar\beta_2,\ldots,\bar\beta_{k}\}
	\end{equation*}
	and proceed with the inductive step assuming that $\mathcal B_{k-1}$ is a basis for $\mathcal G_{K,T_{k-1}}$.  We shall first show that $\mathcal B_k$ is linearly independent.  To this end,
	 we assume that $r_1,r_2,\ldots,r_{n-1},s_1,s_2,\ldots,s_k\in \rat$ are such that
	 \begin{equation} \label{LIStart}
	 	1 = \prod_{i=1}^{n-1}\bar\alpha_i^{r_i} \prod_{j=1}^k \bar\beta_j^{s_j}.
	\end{equation}
	We remind the reader that $\|\cdot\|_{w_k}$ is indeed a well-defined map on $\mathcal G_K$ which satisfies \ref{Multiply} and \ref{Scalar}.  Moreover, by our assumption \ref{NonArchChoices}, we have 
	$\|\beta_j\|_{w_k} = 1$ for all $1\leq j\leq k-1$.  Additionally, the points $\alpha_i$ are ordinary units in $K$, meaning that $\|\alpha_i\|_{w_k} = 1$ for all $1\leq i\leq n-1$.  Now applying $\|\cdot \|_{w_k}$ to both sides
	of \eqref{LIStart}, we conclude that
	\begin{equation*}
		1 = \|\bar\beta_k\|^{s_k}_{w_k}.
	\end{equation*}
	Again by our assumption \ref{NonArchChoices}, we know that $\|\bar\beta_k\|_{w_k} < 1$, and hence, we conclude that $s_k = 0$.  Plugging this value into \eqref{LIStart}, we obtain
	 \begin{equation*} \label{LIInductive}
	 	1 = \prod_{i=1}^{n-1}\bar\alpha_i^{r_i} \prod_{j=1}^{k-1} \bar\beta_j^{s_j}
	\end{equation*}
	so it follows from the inductive hypothesis that $r_i =0$ for all $1\leq i\leq n-1$ and $s_j = 0$ for all $1\leq j\leq k-1$.
	
	Since $\mathcal B_k$ contains $n+k-1$ elements, in order to complete the proof, it is sufficient to show that $\dim(\mathcal G_{K,T_k}) \leq n+k-1$.  To this end, we apply Dirichlet's $S$-unit Theorem with 
	$T_k$ in place of $S$ to obtain a fundamental set of $T_k$-units $\{\gamma_1,\gamma_2,\ldots,\gamma_{n+k-1}\}$.  We claim that
	\begin{equation*}
		\{\bar\gamma_1,\bar\gamma_2,\ldots,\bar\gamma_{n+k-1}\}
	\end{equation*}
	spans $\mathcal G_{K,T_k}$.  To see this, let $\bar\gamma\in\mathcal G_{K,T_k}$ and let $\gamma$ be a representative of $\bar\gamma$ in $K_\divv$.  Then there exists $\ell\in \nat$ such that $\gamma^\ell\in K^\times$,
	and note that $\gamma^\ell$ must belong to $U_{K,T_k}$.  As a result, we obtain integers $r_1,\ldots,r_{n+k-1}$ such that $\gamma^\ell = \prod_{i=1}^{n+k-1} \gamma_i^{r_i}$, and therefore,
	\begin{equation*}
		\bar\gamma = \prod_{i=1}^{n+k-1} \bar\gamma_i^{r_i/\ell}
	\end{equation*}
	proving that $\bar\gamma\in \spann\{\bar\gamma_1,\bar\gamma_2,\ldots,\bar\gamma_{n+k-1}\}$.  It now follows that $\dim(\mathcal G_{K,T_k}) \leq n+k-1$, establishing the lemma.
\end{proof}

We will also need the following basic linear algebra lemma.

\begin{lem} \label{SumsCosets}
	Suppose $V$ is vector space over a field $F$ containing subspaces $A$ and $B$.  Let $x,y\in V$.
	\begin{enumerate}[label=(\roman*)]
		\item\label{Sum} If $A+B = V$ then $(x+A)\cap (y+B)\ne \emptyset$.
		\item\label{Intersect} If $A\cap B = \{0\}$ then $(x+A)\cap (y+B)$ contains at most one point.
		\item\label{DirectSum} If $A\oplus B = V$ then $(x+A)\cap (y+B)$ contains a unique point.
	\end{enumerate}
\end{lem}
\begin{proof}
	For \ref{Sum}, we observe that $x = a_1 + b_1$ and $y = a_2+b_2$ for some $a_1,a_2\in A$ and $b_1,b_2\in B$, and therefore, we obtain
	\begin{equation*}
		x + a_2 = a_1 + a_2 + b_1\quad\mbox{and}\quad y + b_1 = a_2+b_1 + b_2.
	\end{equation*}
	Hence
	\begin{equation*}
		a_2 + b_1 = x + a_2 - a_1 \in x+A\quad\mbox{and}\quad a_2 + b_1 = y+b_1 - b_2 \in y+B,
	\end{equation*}
	showing that $a_2 + b_1 \in (x+A)\cap (y+B)$.
	
	To prove \ref{Intersect}, suppose that $a,b\in (x+A)\cap (y+B)$.  It follows that $a-b$ belongs to both $A$ and $B$, forcing $a-b=0$, or equivalently, $a = b$.  Now \ref{DirectSum} follows from \ref{Sum}
	and \ref{Intersect}.
		
\end{proof}

In order to prove Theorem \ref{SurjectiveOnly}, we will need to identify a consistent map associated to an arbitrary linear map $\Phi:\mathcal G\to \real$.  Our next lemma is the key ingredient in constructing such a map.
For a linear map $\Phi:\mathcal G\to \real$ and $\alpha\in \algt$, we shall adopt the convention that $\Phi(\alpha) = \Phi(\bar\alpha)$.

\begin{lem} \label{AllPhi}
	Suppose $\Phi:\mathcal G\to \real$ is a $\rat$-linear map and $r\in \real$.  Then for every number field $K$, there exists a consistent map $c_K:\mathcal J\to \real$ such that
	\begin{equation} \label{ConsistentChoice}
		c_K(\rat,\infty) = r\quad\mbox{and}\quad \Phi(\alpha) = \Phi_{c_K}(\alpha)\mbox{ for all } \alpha\in K^\times.
	\end{equation}
	Moreover, if $L$ is a finite extension of $K$ and $c_L$ satisfies \eqref{ConsistentChoice} with $L$ in place of $K$, then $c_L(K,v) = c_K(K,v)$ for all places $v$ of $K$.
\end{lem}
\begin{proof}
	We shall select values $y_v$ for each $v\in M_K$ and then apply Theorem \ref{CanonicalConsistent}.  As in the proof of Lemma \ref{ZeroPhi}, we begin by considering the case where $v$
	is Archimedean.  We suppose that $\{v_1,v_2,\ldots,v_n\}$ are the Archimedean places of $K$, assume that $\{\alpha_1,\alpha_2,\ldots,\alpha_{n-1}\}$ is a fundamental set of units in $K$, and let $A$ be the matrix 
	from the statement of Lemma \ref{RankRegulator}.  In this situation, $\ker(A)$ is a one-dimensional subspace of $\real^n$.  Writing 
	\begin{equation*}
		\Lambda = (\lambda(K,v_1),\lambda(K,v_2),\ldots,\lambda(K,v_n))^T,
	\end{equation*}
	we see clearly that $\ker(A) = \{r\Lambda: r\in \real\}$.  We shall further write 
	\begin{equation*}
		\bb = (\Phi(\alpha_1),\Phi(\alpha_2),\ldots,\Phi(\alpha_{n-1}))^T
	\end{equation*}
	and note that $A^{-1}(\bb)$ is a coset of $\ker(A)$ in $\real^n$.
	
	Now define the $n-1$ dimensional subspace $\Delta$ of $\real^n$ by
	\begin{equation*}
		\Delta = \left\{ (x_1,x_2,\ldots, x_n)^T: \sum_{i=1}^n x_i = 0\right\}.
	\end{equation*}
	Plainly we have that $\Delta\cap \ker(A) = \{0\}$ so that $\real^n = \Delta \oplus \ker(A)$.  According to Lemma \ref{SumsCosets}\ref{DirectSum}, each coset of $\Delta$ shares a unique point with each coset of $\ker(A)$.
	Therefore, there exists a unique point $\xx = (x_1,x_2,\ldots,x_n)^T$ such that
	\begin{equation} \label{ArchYDef}
		\xx\in \left[ (r,0,0,\ldots,0) + \Delta\right] \cap A^{-1}(b).
	\end{equation}
	We let $y_{v_i} = x_i$ for all $i$, and as a result, we have defined $y_v$ for each Archimedean place $v$ of $K$ and note that
	\begin{equation}\label{ySum}
		\sum_{v\mid \infty} y_v = r.
	\end{equation}
	
	Now let $\{w_1,w_2,w_3,\ldots\}$ be the complete list of non-Archimedean places of $K$.  For each $j$, we apply \cite[Lemma 3.1]{SamuelsClassification} to obtain a point $\beta_j\in \mathcal O_K$ 
	such that $\|\beta_j\|_{w_j} < 1$ and $\|\beta_j\|_{w_i} = 1$ for all $i\ne j$.  For $m\in \nat\cup\{0\}$, we now define
	\begin{equation*}
		S_m = \left\{v_1,v_2,\ldots,v_n,w_1,w_2,\ldots,w_m\right\}.
	\end{equation*}
	Now for each $j\in \nat$, we define
	\begin{equation} \label{NonArchYDef}
		y_{w_j}  = \frac{\Phi(\beta_{j}) - \sum_{v\mid \infty} y_v \log\|\beta_{j}\|_v}{\log \|\beta_j\|_{w_j}}
	\end{equation}
	so that we have defined $y_v$ for all non-Archimedean places $v$ of $K$.  According to Theorem \ref{CanonicalConsistent}, there exists a consistent map $c_K:\mathcal J\to \real$ such that
	\begin{equation} \label{CDef}
		c_{K}(K,v) = y_v\quad\mbox{for all } v\in M_K.
	\end{equation}
	Because of \eqref{ySum}, $c_K$ clearly satisfies the first equality of \eqref{ConsistentChoice} so it remains only to show the second equality.
	
	To see this, we obtain from \eqref{ArchYDef} that
	\begin{equation*}
		\begin{pmatrix} \Phi(\alpha_1) \\ \Phi(\alpha_2) \\ \vdots \\ \Phi(\alpha_{n-1}) \end{pmatrix} = b = A\xx =   A\begin{pmatrix} y_{v_1} \\ y_{v_2} \\ \vdots \\ y_{v_n} \end{pmatrix}  = 
			\begin{pmatrix} \sum_{i=1}^n y_{v_i} \log\|\alpha_1\|_{v_i} \\ \sum_{i=1}^n y_{v_i} \log\|\alpha_2\|_{v_i} \\ \vdots \\ \sum_{i=1}^n y_{v_i} \log\|\alpha_{n-1}\|_{v_i} \end{pmatrix},
	\end{equation*}
	or equivalently,
	\begin{equation*}
		\begin{pmatrix} \Phi(\alpha_1) \\ \Phi(\alpha_2) \\ \vdots \\ \Phi(\alpha_{n-1}) \end{pmatrix} = 
			\begin{pmatrix} \sum_{v\mid\infty} c_K(K,v) \log\|\alpha_1\|_{v} \\ \sum_{v\mid\infty} c_K(K,v) \log\|\alpha_2\|_{v}  \\ \vdots \\ \sum_{v\mid\infty} c_K(K,v) \log\|\alpha_{n-1}\|_{v}  \end{pmatrix}.
	\end{equation*}
	Each point $\alpha_i$ is a unit, so we have $\|\alpha_i\|_v = 1$ for all non-Archimedean places $v$ of $K$. Therefore, it follows that
	\begin{equation*}
		\begin{pmatrix} \Phi(\alpha_1) \\ \Phi(\alpha_2) \\ \vdots \\ \Phi(\alpha_{n-1}) \end{pmatrix} = \begin{pmatrix} \Phi_{c_K}(\alpha_1) \\ \Phi_{c_K}(\alpha_2) \\ \vdots \\ \Phi_{c_K}(\alpha_{n-1}) \end{pmatrix},
	\end{equation*}
	or in other words,
	\begin{equation} \label{ArchMatch}
		\Phi(\alpha_i) = \Phi_{c_K}(\alpha_i)\quad\mbox{for all } 1\leq i\leq n-1.
	\end{equation}
	
	Next, we let $j\in \nat$ and observe that
	\begin{align*}
		\Phi_{c_K}(\beta_j) & = \sum_{v\in M_K} c(K,v) \log \|\beta_j\|_v \\
						& = \sum_{i=1}^\infty c(K,w_i) \log \|\beta_j\|_{w_i} + \sum_{v\mid\infty} c(K,v) \log \|\beta_j\|_v \\
						& = c(K,w_j) \log\|\beta_j\|_{w_j} + \sum_{v\mid\infty} c(K,v) \log \|\beta_j\|_v,
	\end{align*}
	where the last equality follows because of our assumption that $\|\beta_j\|_{w_i} = 1$ for all $i\ne j$.  Now using \eqref{CDef}, we obtain
	\begin{equation*}
		\Phi_{c_K}(\beta_j) = y_{w_j}\log \|\beta_j\|_{w_j} + \sum_{v\mid\infty} y_v \log \|\beta_j\|_v,
	\end{equation*}
	and using \eqref{NonArchYDef} to replace $y_{w_j}$, we quickly obtain
	\begin{equation}\label{NonArchMatch}
		\Phi(\beta_j) = \Phi_{c_K}(\beta_j)\quad\mbox{for all } j\in \nat.
	\end{equation}
	
	Assuming now that $\alpha\in K^\times$, there exists $m\in \nat$ such that $\alpha$ is an $S_m$-unit. In other words, we have $\bar\alpha\in \mathcal G_{K,S_m}$.  According to Lemma \ref{Dirichlet},
	$\{\bar\alpha_1,\bar\alpha_2,\ldots,\bar\alpha_{n-1},\bar\beta_1,\bar\beta_2,\ldots,\bar\beta_m\}$ is a basis for $\mathcal G_{K,S_m}$, so there exist $r_1,\ldots,r_{n-1},s_1,\ldots s_m\in\rat$ such that
	\begin{equation*}
		\bar\alpha = \prod_{i=1}^{n-1}\bar\alpha_i^{r_i} \prod_{i=1}^m \bar\beta_j^{s_j}.
	\end{equation*}
	Since $\Phi$ and $\Phi_{c_K}$ are both linear maps, we conclude that $\Phi(\bar\alpha) = \Phi_{c_K}(\bar\alpha)$, as required.
	
	For the second statement of the Lemma, we have assumed that $\Phi_{c_L}(\alpha) = \Phi(\alpha)$ for all $\alpha \in L^\times$, so in particular, this equality holds for all $\alpha\in K^\times$.
	As a result, we have
	\begin{equation*}
		\Phi_{c_L}(\alpha) = \Phi_{c_K}(\alpha) \quad\mbox{for all } \alpha\in K^\times.
	\end{equation*}
	For each such point $\alpha$, Theorem \ref{KernelOnly} implies that $0 = \Phi_{c_L}(\alpha) - \Phi_{c_K}(\alpha) = \Phi_{c_L - c_K}(\alpha)$.  If $v$ is a place of $K$, we may apply Lemma \ref{ZeroPhi}
	with $c = c_K - c_L$ to conclude that
	\begin{equation*}
		[c_L - c_K](K,v) = [c_L(\rat,\infty) - c_K(\rat,\infty)] \lambda(K,v) = (r-r)\lambda(K,v) = 0,
	\end{equation*}
	so the result follows immediately.
\end{proof}

We now complete the proof of Theorem \ref{RealSurjective} by proving that $\Phi^*$ is surjective.

\begin{proof}[Proof of Theorem \ref{SurjectiveOnly}]
	We suppose that $\Phi:\mathcal G\to \real$ is a $\rat$-linear map.  For each number field $K$, we apply Lemma \ref{AllPhi} with $r=0$ to obtain a consistent map $c_K:\mathcal J\to \real$ such that
	\begin{equation} \label{PartialPhi}
		c_K(\rat,\infty) = 0\quad\mbox{and}\quad \Phi(\alpha) = \Phi_{c_K}(\alpha)\mbox{ for all } \alpha\in K^\times.
	\end{equation}
	Now we let $c:\mathcal J\to \real$ be given by
	\begin{equation*}
		c(K,v) = c_K(K,v)
	\end{equation*}
	and claim that $c$ is consistent.  To this end, we assume that $K$ is a number field, $v$ is a place of $K$, and $L$ is a finite extension of $K$.  Using the fact that $c_L$ is consistent, we obtain that
	\begin{equation*}
		\sum_{w\mid v} c(L,w) = \sum_{w\mid v} c_L(L,w) = c_L(K,v),
	\end{equation*}
	and then applying the second statement of Lemma \ref{AllPhi}, we conclude that
	\begin{equation*}
		\sum_{w\mid v} c(L,w) = c_K(K,v) = c(K,v).
	\end{equation*}
	Since $c$ is now known to be consistent, the function $\Phi_c$ is well-defined.  Moreover, if $\alpha\in \algt$, we let $K$ be a number field containing $\alpha$ and apply the first statement of Lemma \ref{AllPhi}
	to obtain
	\begin{equation*}
		\Phi_c(\alpha) = \sum_{v\in M_K} c(K,v) \log\|\alpha\|_v = \sum_{v\in M_K} c_K(K,v) \log\|\alpha\|_v = \Phi_{c_K}(\alpha) = \Phi(\alpha).
	\end{equation*}
	As this equality holds for all $\alpha\in \algt$, we have shown that $\Phi^*(c) = \Phi$ implying that $\Phi^*$ is surjective.	
\end{proof}

Our only remaining task is to prove Theorem \ref{KRational}.

\begin{proof}[Proof of Theorem \ref{KRational}]
	 For simplicity, we let $\Phi:K^\times\to \real$ be given by
	\begin{equation*}
		\Phi(\alpha) = \sum_{v\in M_K} y_v\log\|\alpha\|_v
	\end{equation*}
	and first assume that $\Phi(\alpha)\in\rat$ for all $\alpha\in K^\times$.  Each point $\alpha_i$ is a unit in $K$, and hence, $\|\alpha_i\|_w = 1$ for all non-Archimedean places $w$ of $K$.  For all $1\leq i\leq n$,
	we now have that
	\begin{equation*}
		\sum_{i=1}^n y_i \log\|\alpha_i\|_{v_i} = \Phi(\alpha_i) \in \rat,
	\end{equation*}
	and \ref{Arches} follows immediately.  Similarly, we have assumed that $\|\beta_v\|_w = 1$ for all non-Archimedean places $w\ne v$, so we conclude immediately that \ref{NonArches} holds.
	
	Now suppose that \ref{Arches} and \ref{NonArches} hold and let $\alpha\in K^\times$.  Let $S_\infty$ be the set of Archimedean places of $K$.
	There exists a finite set $S$ of places of $K$ containing $S_\infty$ such that $\alpha\in U_{K,S}$, and hence, $\bar\alpha\in \mathcal G_{K,S}$.  Let $S_0 = S\setminus S_\infty$.
	Now according to Lemma \ref{Dirichlet}, the set
	\begin{equation*}
		\left\{\bar\alpha_1,\bar\alpha_2,\ldots,\bar\alpha_{n-1}\right\} \cup \left\{\beta_v:v\in S_0\right\}
	\end{equation*}
	is a basis for $\mathcal G_{K,S}$.  Therefore, there exist $r_1,r_2,\ldots,r_{n-1}\in \rat$ and $s_v\in \rat$ for all $v\in S_0$ such that
	\begin{equation*}
		\bar\alpha = \prod_{i=1}^{n-1}\bar\alpha_i^{r_i} \prod_{v\in S_0} \bar\beta_v^{s_v}.
	\end{equation*}
	Now using our definitions presented prior to Lemma \ref{Dirichlet}, $\Phi$ is a well-defined linear map on $\mathcal G_{K,S}$, and hence
	\begin{equation*}
		\Phi(\alpha) = \Phi(\bar\alpha) = \sum_{i=1}^{n-1} r_i\Phi(\bar\alpha_i) + \sum_{v\in S_0} s_v\Phi(\bar\beta_v) = \sum_{i=1}^{n-1} r_i\Phi(\alpha_i) + \sum_{v\in S_0} s_v\Phi(\beta_v)
	\end{equation*}
	It follows from \ref{Arches} that $\Phi(\alpha_i) \in \rat$ for all $i$, and \ref{NonArches} implies that $\Phi(\beta_v) \in \rat$ for all $v\in S_0$.  We now obtain $\Phi(\alpha)\in \rat$, as required.
\end{proof}

\section{Acknowledgement}

We thank the referee for many useful comments regarding this work, especially several important corrections to the proof of Theorem \ref{CanonicalConsistent}.  We additionally thank
Christopher Newport University for the Spring 2025 sabbatical leave during which the majority of this work took place.  Finally, we recognize the AMS-Simons Research Enhancement Grant for PUI 
Faculty for partially funding this work.

\bibliographystyle{abbrv}
\bibliography{GRepresentations}

\end{document}